\documentclass[12pt]{amsart}
\usepackage{amsmath,amscd,amssymb,amsfonts}
\theoremstyle{plain}
\newtheorem{theorem}{Theorem}[section]
\newtheorem{proposition}[theorem]{Proposition}

\theoremstyle{definition}

\newtheorem{remark}[theorem]{Remark}
\newtheorem{example}[theorem]{Example}

\newcommand{\msum}{\hbox{$\sum$}}
\newcommand{\mprod}{\hbox{$\prod$}}
\newcommand{\bA}{{\mathbb A}}
\newcommand{\bC}{{\mathbb C}}
\newcommand{\bN}{{\mathbb N}}
\newcommand{\bQ}{{\mathbb Q}}
\newcommand{\bR}{{\mathbb R}}
\newcommand{\bZ}{{\mathbb Z}}
\newcommand{\bc}{{\mathbf c}}
\newcommand{\bs}{{\mathbf s}}

\newcommand{\cD}{{\mathcal D}}
\newcommand{\cO}{{\mathcal O}}
\newcommand{\fa}{\mathfrak{a}}

\begin{document}

\title[Roots of Bernstein-Sato polynomials]
{Combinatorial description of the roots\\
of the Bernstein-Sato polynomials\\
for monomial ideals}

\author[N. Budur]{Nero Budur}
\address{Department of Mathematics, The Johns Hopkins
University, Baltimore, MD 21218, USA}
\email{{\tt budur@math.jhu.edu}}

\author[M. Musta\c{t}\v{a}]{Mircea~Musta\c{t}\v{a}}
\address{Department of Mathematics, University of Michigan,
Ann Arbor, MI 48109, USA}
\email{{\tt mmustata@umich.edu}}

\author[M. Saito]{Morihiko Saito}
\address{RIMS Kyoto University, Kyoto 606-8502, Japan}
\email{{\tt msaito@kurims.kyoto-u.ac.jp}}

\begin{abstract}
We give a combinatorial description of the roots of the
Bernstein-Sato polynomial of a monomial ideal using
the Newton polyhedron and some semigroups
associated to the ideal.
\end{abstract}

\thanks{The research of the second author was partially supported by
NSF grant 0500127}

\maketitle

\section{Introduction}

The Bernstein-Sato polynomial 
of a function was introduced
independently by J.~Bernstein and M.~Sato.
It was generalized in \cite{BMS2} to the case of arbitrary ideals in a
polynomial ring (see also \cite{Gyo} and \cite{Sab} for the case where the
generators of the ideal are chosen).
For monomial ideals, it was shown in \cite{BMS2} that in principle, this
polynomial can be computed algorithmically. 
However, an explicit combinatorial description of its roots
was missing.
In this note we give such a combinatorial description of the roots
(without multiplicities)
of the Bernstein-Sato polynomial of a
monomial ideal, using the Newton polyhedron and some semigroups
associated to the ideal.
This description was first obtained as a consequence of the 
main theorems in \cite{BMS1} that used reduction mod $p$.
In this paper we give a direct proof.

In order to state our first main theorem, we need to introduce some notation.
Let
$ \fa $ be a nonzero monomial ideal in the polynomial ring
$ \bC[x] = \bC[x_{1},\dots, x_{n}] $. We denote the monomial corresponding
to $u\in\bZ_{\geq 0}^n$ by $x^u$.
Let $\Gamma_{\fa} \subset \bZ^{n}_{\geq 0}$ be the semigroup corresponding to
$ \fa $, i.e.
$ \Gamma_{\fa} = \{u\in\bN^n\mid x^u\in\fa\}$.
The Newton polyhedron
$ P_{\fa} $ of
$ \fa $ is the convex hull of $\Gamma_{\fa}$ in
$ \bR_{\ge 0}^{n} $.

Let $e=(1,\ldots,1)$ in $\bZ_{\geq 0}^n$.
If $Q$ is any face of $P_{\fa}$,
we denote by $M_{Q}$ the translate by $e$ of the subsemigroup of $\bZ^n$
generated by elements of the form $u-v$, where $u$ is in $\Gamma_{\fa}$
and $v$ is in $\Gamma_{\fa}\cap Q$.
If $v_0$ is an element in $\Gamma_{\fa}\cap Q$, then
we put $M_Q':=v_0+M_Q$ (this is a subset 
of $M_Q$ that does not depend
on the choice of $v_0$).

If a face $Q$ of $P_{\fa}$ is not contained in any coordinate
hyperplane, then its affine span does not contain the
origin (see Remark~\ref{rem5}),
and there is a (not necessarily unique) linear function $L_Q$ on
$ \bR^{n} $ with rational coefficients
such that
$L_Q=1$ on $Q$. 
Let $V_Q$ be the linear subspace
generated by $Q$.
Put
$$R_Q:=\{-L_Q(u)\mid u\in (M_Q\smallsetminus M_{Q}')\cap V_Q\}.$$
Note that
$ (M_Q\smallsetminus M_{Q}')\cap V_Q $ is empty if
$ Q $ is contained in a coordinate hyperplane.
If $Q$ is a facet (i.e. a maximal-dimensional proper face) of $P_{\fa}$, then $L_Q$ is
unique
and we denote by
$ m_Q $ the smallest positive integer such that
$ m_QL_Q $ has integral coefficients.

With this notation, we prove the following

\begin{theorem}\label{main1}
The roots of the Bernstein-Sato polynomial of $\fa$ are the union of
the sets $ R_Q $ for the faces $Q$ of $P_{\fa}$ that are not contained
in any coordinate hyperplanes.
\end{theorem}

If $n\geq 2$ it is enough to restrict to faces
of dimension at least one. On the other hand, it is not enough to
consider only the facets of $P_{\fa}$. 

As a consequence of the above theorem, we obtain the following description 
of the roots of the Bernstein-Sato polynomial of $\fa$, modulo $\bZ$.
Note that in this case the description depends only on the Newton polyhedron
of the ideal, i.e. only on the integral closure of $\fa$.
See \cite{BMS1} for an approach to this theorem via characteristic $p$ methods.

\begin{theorem}\label{main2}
The set of classes in $\bQ/\bZ$
 of the roots of the Bernstein-Sato polynomial of $\fa$
 is equal to the union of the subgroups generated by $1/m_Q$,
with $Q$ running over the facets of $P_{\fa}$ that are not
contained in coordinate hyperplanes.
\end{theorem}

\medskip

In \S 2 we recall the definition of Bernstein-Sato polynomials,
with emphasis on the case of monomial ideals.
In \S 3 we prove the main theorems, and
in the last section we give some examples to
illustrate our combinatorial description.

\section{Bernstein-Sato polynomials}

In this section, we recall the definition of Bernstein-Sato
polynomials and interpret this definition in the special case
of a monomial ideal. For details we refer to \cite{BMS2}.

Let
$ X $ be a smooth affine variety over
$ \bC $, and let
$ Z $ be a (not necessarily reduced or irreducible) subvariety of
$ X $, different from $X$.
Fix generators $f_1,\ldots,f_r$ for the ideal of
$ Z $.
We denote by
$ \cD_{X} $ the sheaf of linear differential operators on
$ X $.
It acts naturally on
$$
\cO_{X}[\mprod_{i}f_{i}^{-1},s_{1},\dots, s_{r}]
\mprod_{i}f_{i}^{{s}_{i}},
$$
where the
$ s_{i} $ are independent variables.
We define a
$ \cD_{X} $-linear action of
$ t_{i} $ by
$ t_{i}(s_{j}) = s_{j} + 1 $ if
$ j = i $, and
$ t_{i}(s_{j}) = s_{j} $ otherwise.
Let
$ s_{i,j} = s_{i}t_{i}^{-1}t_{j} $
and
$ s = \msum_{i}s_{i} $.
The Bernstein-Sato polynomial (or the
$ b $-function)
$ b_{f}(s) $ of
$ f := (f_{1}, \dots, f_{r}) $ is defined to be the monic
polynomial of the lowest degree in
$ s $ satisfying a relation of the form
\begin{equation}\label{eq1}
b_{f}(s)\mprod_{i}f_{i}^{{s}_{i}} =
\msum_{k=1}^r P_{k}t_{k}
\mprod_{i}f_{i}^{{s}_{i}},
\end{equation}
with the
$ P_{k} $ in the ring generated by
$ \cD_{X} $ and the
$ s_{i,j} $.

\begin{remark}\label{rem1}
We can identify
$ f_{i}^{s_{i}} $ with the delta function
$ \delta(t_{i}-f_{i}) $, so that
$ s_{i}t_{i}^{-1} $ corresponds to
$ -\partial_{t_{i}} $, see \cite{Mal}.
Therefore
$ \msum_{k=1}^r P_{k}t_{k} $ in the above definition can be
replaced with an element of the ring generated by
$ \cD_{X} $ and
$ \prod_{i} t_{i}^{\mu_{i}}\partial_{t_{i}}^{\nu_{i}} $ with
$ \sum_{i}\mu_{i} - \sum_{i}\nu_{i} \ge 1 $.
\end{remark}

\begin{remark}\label{rem2}
Using the theory of
$ V $-filtrations of Kashiwara and Malgrange, we can show that
$ b_{f}(s) $ is independent of the choice of the generators of
the ideal and depends only on the variety
$ Z $, see \cite{BMS2}.
Furthermore, the roots of
$ b_{f}(s) $ are negative rational numbers, generalizing \cite{Kas}.
Note that
$ b_{f}(s) $ coincides with the polynomial
$ b_{\alpha} $ for the index $\alpha=(1,\ldots,1)$,
that appeared in \cite{Gyo}, 2.13
(see also \cite{Sab}, I, 3.1).
However,
$ b_{f}(s) $ is slightly different from the
 polynomial in \cite{Sab}, II,
Prop.~1.1,
because our definition requires certain additional binomial
polynomials as in (\ref{eq2}) below.
We mention also that Sabbah proved in \cite{Sab}
the existence of nonzero polynomials of
{\it several variables}
that satisfy functional equations similar to (\ref{eq1}) above.
\end{remark}

\medskip

We give now an equivalent definition of the Bernstein-Sato polynomial.
For
$ \bc = (c_{1}, \dots, c_{r}) \in \bZ^{r} $, let
$ I(\bc)_{-} = \{i\mid c_{i} < 0 \} $.
The Bernstein-Sato polynomial
$ b_{f}(s) $ is the monic polynomial of the smallest degree
such that
$ b_{f}(s)\mprod_{i}f^{s_{i}} $ belongs to the
$ \cD_{X}[s_{1}, \dots, s_{r}] $-submodule generated by
\begin{equation}\label{eq2}
\mprod_{i\in I(\bc)_{-}}\hbox{$\binom{s_{i}}{-c_{i}}$}\cdot
\mprod_{i=1}^{r}f_{i}^{s_{i}+c_{i}},
\end{equation}
where
$ \bc = (c_{1}, \dots, c_{r}) $ runs over the elements of
$ \bZ^{r} $ such that
$ \msum_{i}c_{i} = 1 $.
Here
$ s = \msum_{i=1}^{r}s_{i} $ and
$ \binom{s_{i}}{m} = s_{i}(s_{i}-1)\cdots (s_{i}-m+1)/m! $.

This definition of the Bernstein-Sato polynomial coincides with
the previous one.
Indeed, we have the relation
$ t_{i}^{-1}s_{i} = (s_{i}-1)t_{i}^{-1} $, which implies
$$
(s_{i}t_{i}^{-1})^{-c_{i}} = (-c_{i})!\hbox{$\binom{s_{i}}{-c_{i}}$}
t^{c_{i}} \quad\text{for}\quad c_{i} < 0.
$$
Letting
$ \theta_{i} = t_{i} $ if
$ c_{i} > 0 $, and
$ \theta_{i} = \partial_{t_{i}}^{-1} $ if
$ c_{i} < 0 $, it is enough to consider
$ \prod_{i}\theta_{i}^{c_{i}} $ for
$ \bc \in \bZ^{r} $ with
$ \sum_{i}c_{i} = 1 $.

In the case of a monomial ideal we can make
more explicit the above definition.
Assume that
$ X $ is the affine space
$ \bA^{n} $ and the
$ f_{j} $ are monomials with respect to the coordinate
system
$(x_{1},\dots,x_{n}) $ on
$ \bA^{n} $.
Write
$ f_{j} = \prod_{i=1}^{n}x_{i}^{a_{i,j}} $ and let
$\ell_{i}(\bs) = \msum_{j=1}^{r}a_{i,j}s_{j}$.
for
$ \bs = (s_{1},\dots,s_{r}) $. Therefore we have
\begin{equation}\label{eq3}
\mprod_{j=1}^{r}f_{j}^{s_{j}} =
\mprod_{i=1}^{n}x_{i}^{\ell_{i}(\bs)}.
\end{equation}
We put
$ \ell(\bc) = (\ell_{1}(\bc),\dots,\ell_{n}(\bc)) $ and
$ I'(\ell(\bc))_{+} = \{i\mid\ell_{i}(\bc) > 0\} $.
For every $\bc$ in $\bZ^r$ such that $\msum_jc_j=1$, we define
\begin{equation}\label{eq4}
g_{\bc}(s_{1},\dots,s_{r}) =
\prod_{j\in I(\bc)_{-}}\binom{s_{j}}{-c_{j}}\cdot
\prod_{i\in I'(\ell(\bc))_{+}}
\binom{\ell_{i}(\bs)+\ell_{i}(\bc)}{\ell_{i}(\bc)}.
\end{equation}
Let
$ I_{\fa} $ be the ideal of
$\bQ[s_{1},\dots,s_{r}] $ generated by
$ g_{\bc}(s_{1},\dots,s_{r}) $ where
$ \bc = (c_{1},\dots,c_{r}) $ runs over the elements of
$ \bZ^{r} $ with
$ \msum_{i}c_{i} = 1 $.
Using the
$ \bZ^{n} $-grading on
$ \cD_{X} $ such that the degree of
$ x_{i} $ is the
$ i $th unit vector
$ e_{i} $ of
$ \bZ^{n} $, and the degree of
$ \partial_{x_{i}} $ is
$ -e_{i} $, one can show the following

\begin{proposition}\label{prop1}{\rm (}\cite{BMS2}{\rm )}
With the above notation,
the Bernstein-Sato polynomial
$ b_{f}(s) $ is the monic polynomial of smallest degree
such that
$ b_{f}(\sum_{i}s_{i}) $ belongs to the ideal
$ I_{\fa} $.
\end{proposition}

\begin{remark}\label{rem3}
It is not necessarily easy to give finite generators of the ideal
$ I_{\fa} $ explicitly.
One problem is to give a bound for the  corresponding $\bc$
in terms of the $a_{i,j}$. An algorithm for finding finite
generators was given
in \cite{BMS2}, but it is not an easy task to write down
an explicit program even if
$ n $ and
$ r $ are small.
Once a finite system of
 generators for the ideal is given, it is not very
difficult to determine all the roots 
of the Bernstein-Sato polynomial (without multiplicities).
\end{remark}

\section{Proofs of the main theorems}

In this section, we prove Theorems~\ref{main1} and \ref{main2}.
Before giving the proof, we recall the notation introduced in the 
Introduction and make some preliminary remarks.

Let $\fa\subseteq\bC[x_1,\ldots,x_n]$ be an ideal generated by
$f_1,\ldots,f_r$, where $f_j=x^{v_j}$
and $v_j=(a_{1,j},\ldots,a_{n,j})$.
We denote by $\Gamma_{\fa}$ the set of those $u$ in $\bN^n$
such that $x^u$ is in $\fa$.
Consider a proper face $Q$ of the Newton polyhedron $P_{\fa}$
of $\fa$.
 
Recall that if $e=(1,\ldots,1)\in\bZ^n$, then $M_{Q}$ is defined such that
$M_Q-e$ is the subsemigroup of $\bZ^n$ generated by
$$\{u-v\mid u\in \Gamma_{\fa}, v\in\Gamma_{\fa}\cap Q\}.$$
We have also defined $M_Q':=v_0+M_Q$, where $v_0$ is an arbitrary element in
$\Gamma_{\fa}\cap Q$. More generally, for
$k$ in $\bZ_{\geq 0}$ we put $M_Q^{(k)}:=kv_0+M_Q$. It is clear that the definition
does not depend on the choice of $v_0$. Moreover, since
$\Gamma_{\fa}+\bZ^n_{\geq 0}\subseteq\Gamma_{\fa}$, we deduce
$M_Q^{(k)}+\bZ_{\geq 0}^n\subseteq M_Q^{(k)}$ for every $k$.

Suppose now that the affine span of $Q$ does not contain the origin,
so we have a linear function $L_Q$ on $\bR^n$ having $\bQ$-coefficients such that
$L_Q=1$ on $Q$. Note that $R_Q$ is computed by looking at values of
$L_Q$ on the linear subspace $V_Q$ generated by $Q$, so it is independent
of our choice of $L_Q$. 

\begin{remark}\label{rem6}
It follows from the definition that
for every $k\in\bZ_{\geq 0}$ we have
$$R_Q=\{-L_Q(u)+k\mid u\in (M_Q^{(k)}\smallsetminus M_Q^{(k+1)})\cap V_Q\}.$$
\end{remark}

A key step in the proof of Theorem~\ref{main1} 
is based on induction on $n$, reducing
the statement about the monomial ideal $\fa$ to that for a monomial ideal 
$\fa'$ in fewer variables, such that $P_{\fa'}$ is a suitable projection of
$P_{\fa}$. The next remark deals with the combinatorial aspect of
this reduction.

\begin{remark}\label{rem4}
If $Q$ is an unbounded face of $P_{\fa}$, then there is $i$
such that $Q+e_i\subseteq Q$, where $e_i$ is the $i$th vector of the
standard basis of $\bZ^n$. After renumbering the coordinates, we may assume
that $i=n$. Consider the ideal $\fa'$ in $\bC[x_1,\ldots,x_{n-1}]$
defined by $f'_1,\ldots,f'_r$, where $f'_i=\prod_{i=1}^{n-1}x_i^{a_{i,j}}$.
If $p\colon\bR^n\to\bR^{n-1}$ is the projection onto the first 
$(n-1)$ coordinates, we see that $P_{\fa'}=p(P_{\fa})$. 
Since $Q+e_n\subseteq Q$, one can check that
$Q':=p(Q)$ is a face of $P_{\fa'}$ such that $Q=P_{\fa}\cap p^{-1}(Q')$.
Moreover, this construction gives a bijection between
the faces of $P_{\fa'}$ and the faces $Q$ of $P_{\fa}$
such that $Q+e_n\subseteq Q$.

Note that $Q$ is contained in a coordinate
hyperplane if and only if $Q'$ has the same property.
We also see that $0$ lies in the affine span of $Q$ if
and only if the same holds for $Q'$.
If $0$ is not in the affine span of $Q$, then we can choose
our linear functions such that $L_{Q}=L_{Q'}\circ p$. 
Moreover, it is an easy exercise to show, using the definitions,
that $p(M_Q\smallsetminus M_Q')=M_{Q'}\smallsetminus M_{Q'}'$
and $V_Q=p^{-1}(V_{Q'})$. Therefore we have $R_Q=R_{Q'}$.
\end{remark}

\begin{remark}\label{rem5}
Let $Q$ be an arbitrary face of $P_{\fa}$. Since 
$Q$ is the intersection of its affine span with $P_{\fa}$, we see that
if $0$ lies in this affine span, then for any
$ \alpha \ge 1 $ and
$ u\in Q $ we have
$ \alpha u \in Q $.
We deduce that if $Q$ is bounded and $Q\neq\{0\}$ (note that
$Q=\{0\}$ implies $\fa=\bC[x_1,\ldots,x_n]$), then 
$0$ does not lie in the affine span of $Q$.

Using repeatedly the construction in Remark~\ref{rem4}, we deduce
that if the affine span of a face $Q$ of $P_{\fa}$ contains $0$
then $Q$ is contained in a coordinate hyperplane.
\end{remark}

We can give now the proof of our first main theorem.

\begin{proof}[Proof of Theorem~\ref{main1}]
We proceed by induction on
the number $n$ of variables. We divide the proof into several steps.

\smallskip

\noindent{\emph{Step 1}}. We show that the statement of the theorem holds
for $n=1$. If $\fa=\bC[x]$, then $b_{\fa}(s)=1$, so there are no roots.
Since in
this case $P_{\fa}=\bR_{\geq 0}$, there is only one proper face
$ 0 $, which is contained in a coordinate hyperplane,
so the theorem is satisfied. If $b_{\fa}$ is generated by $x^m$
for some $m\geq 1$, then 
$$b_{\fa}(s)=\prod_{i=1}^m\left(s+\frac{i}{m}\right).$$
On the other hand, $P_{\fa}=\{u\in\bR\mid u\geq m\}$.
The only face we have to consider is $Q=\{m\}$ and we can take $L_Q(u)=u/m$.
Moreover, we have $M_Q=\bZ_{>0}$ and $M_Q'=\bZ_{>m}$, so our
statement follows.

\smallskip

\noindent{\emph{Step 2}}. We use the description 
of the Bernstein-Sato polynomial
of a monomial ideal from \S 2 to give an interpretation of the roots that
we will use from now on. Recall that for $\bs=(s_1,\ldots,s_r)$
we put $\ell_i(\bs)=\sum_{i=1}^na_{i,j}s_j$.

It follows from Proposition~\ref{prop1} that
in order to determine the roots of $b_{\fa}$,
it is enough to determine those $\bs$ such that 
$-\bs\in V(I_{\fa})$. More precisely,
we need to consider $\bs\in\bQ^r$ with the following property:
for every $\bc=(c_1,\ldots,c_r)\in\bZ^r$ with $\sum_jc_j=1$,
either there is $j$ such that 
$$c_j<0\,\,{\rm and}\,\,s_j\in\{c_j+1,\ldots,0\},
\leqno(A)$$
or there is $i$ such that 
$$\ell_i(\bc)>0\,\,{\rm and}\,\,\ell_i(\bs)\in
\{1,\ldots,\ell_i(\bc)\}.
\leqno(B)$$
For every such $\bs$, we have the root 
$-\sum_js_j$ of $b_{\fa}$. Moreover,
every root of $b_{\fa}$ appears in this way.
We study now in more detail the above condition for $\bs$ so that
$-\bs$ is in $V(I_{\fa})$.

\noindent{\emph{Step 3}}. 
First, we use the induction hypothesis to show that it is enough
to consider only those
$\bs\in\bQ^r$ such that $\ell_i(\bs)\in\bZ_{>0}$ for
every $i$. Indeed,
let $\fa'\subseteq\bC[x_1,\ldots,x_n]$ be defined as
in Remark~\ref{rem4}, so for $i\leq n-1$ the linear
functions $\ell_i$ are the same for both ideals. 
Assume
$\ell_n(\bs)\not\in\bZ_{>0}$ so that
$(B)$ for
$i=n$ does not hold for any $\bc$.
Then $-\bs$ is in
$V(I_{\fa'})$ if and only if $-\bs$ is in $V(I_{\fa})$.
Moreover, the
corresponding roots of $b_{\fa'}$ and $b_{\fa}$ are equal. 

The inductive hypothesis together
with Remark~\ref{rem4} show that from now on we may restrict our
attention to those $-\bs$ in $V(I_{\fa})$ such that $\ell_i(\bs)\in
\bZ_{>0}$ for all $i$. Equivalently, we may assume that 
$\sum_{j=1}^rs_jv_j$ is in $\bZ_{>0}^n$.

\smallskip

\noindent{\emph{Step 4}}.
Given $\bs\in\bQ^r$ such that $\sum_js_jv_j$ is in $\bZ_{>0}^n$,
we want to reinterpret the condition in \emph{Step 2}
that $-\bs$ is in $V(I_{\fa})$. 
We write $\bs=\bs'-\bs''$, where
$$
s'_j = 
\begin{cases}
s_j &\text{if $s_j\not\in\bZ_{\leq 0}$,} \\
0 &\text{otherwise.}
\end{cases}
$$
Note that $s''_{j}$ is in $\bZ_{\geq 0}$ for every $j$.
In particular, $k(\bs):=\sum_js''_j$ is in $\bZ_{\geq 0}$.
We consider also the subset of $\{1,\ldots,r\}$
given by $J'(\bs)=\{j\vert s_j\not\in\bZ_{\leq 0}\}$
and its complement $J''(\bs)$.

The condition in \emph{Step 2} says that
$-\bs$ is in $V(I_{\fa})$ if and only if for every
$\bc\in\bZ^r$ with $\sum_jc_j=1$ and $c_j\geq s_j$
for every $j$ in $J''(\bs)$, there is $i$ such that $\ell_i(\bs)
\leq\ell_i(\bc)$.
After rewriting this condition for $\bc':=\bc+\bs''$ instead of $\bc$,
we get the following: $-\bs$ is in $V(I_{\fa})$ if and only if
for every $\bc'$ in $\bZ^r$ with $\sum_jc'_j=k(\bs)+1$
and $c'_j\geq 0$ for $j\in J''(\bs)$, there is $i$ such that
$\ell_i(\bs')\leq\ell_i(\bc')$.

For an arbitrary subset $J$ of $\{1,\ldots,r\}$ and for $k\in\bZ_{\geq 0}$
we denote by
 $W(J,k)$ the set
$$\left\{u\in \bZ^n\vert\,{\rm there}\,{\rm is}\,\bc'\in\bZ^r\,{\rm with}
\,\sum_jc'_j=k,\,c'_j\geq 0\,{\rm for}\,j\in J, u-\sum_jc'_jv_j
\in\bZ_{>0}^n\right\}.$$
With this notation, we see that
$-\bs$ is in $V(I_{\fa})$ if and only if $\sum_js'_jv_j$ is not in
$W(J''(\bs),k(\bs)+1)$.

Note that in any case, $\sum_js'_jv_j$ lies in $W(J''(\bs),k(\bs))$.
Indeed, running the above argument with $k(\bs)$ instead of $k(\bs)+1$,
we see that it is enough to show that there is $\bc\in\bZ^r$ such that
$\sum_jc_j=0$ and $c_j\geq s_j$ for all $j\in J''(\bs)$,
and $\ell_i(\bs)-\ell_i(\bc)\in\bZ_{>0}$ for all $i$. By our assumption
on $\bs$, we may take $\bc=0$.

\smallskip

\noindent{\emph{Step 5}}.
Given $\bs\in\bQ^r$ such that
$\sum_js_jv_j\in \bZ^n_{>0}$, let $Q$ be the smallest face of $P_{\fa}$
containing all $v_j$ with $j\in J'(\bs)$. We show that if $-\bs$ is in
$V(I_{\fa})$, then $Q$ is a proper face of $P_{\fa}$. 
Indeed, if $Q=P_{\fa}$, then we can find $\alpha_j\in\bQ_{\geq 0}$
for $j\in J'(\bs)$ with $\sum_j\alpha_j=1$ such that
$v:=\sum_{j\in J'(\bs)}\alpha_jv_j$ is in the interior of $P_{\fa}$. 
In this case, there is $u$ in the convex hull of all the $v_j$, such that
$v-u$ is in $\bQ_{>0}^n$.
We write $u=\sum_{j=1}^r\beta_jv_j$ for $\beta_j\in\bQ_{\geq 0}$
with $\sum_j\beta_j=1$. Fix an  element $j_0$ in
$J'(\bs)$, take a positive integer $m$ that
is large and divisible enough and consider for every $k\in\bZ_{\geq 0}$
$$kv_{j_0}+m(u-v)=kv_{j_0}+\sum_{j=1}^rm\beta_jv_j
-\sum_{j\in J'(\bs)}m\alpha_jv_j.$$
When $m$ goes to infinity, the entries of this element are negative and
go to infinity in absolute value, so we deduce that
$W(J''(\bs),k)=\bZ^n$ for every $k$. If $k=k(\bs)+1$, this contradicts the
fact that $-\bs\in V(I_{\fa})$, by \emph{Step 4}.

\smallskip

\noindent{\emph{Step 6}}.
Suppose that $\bs\in\bQ^r$ is such that $\sum_js_jv_j$ is in 
$\bZ_{>0}^n$, and let $Q$ be the smallest face of $P_{\fa}$
containing all $v_j$ with $j\in J'(\bs)$. We show that for every
$k\in\bZ_{\geq 0}$, we have $W(J''(\bs),k)=M_Q^{(k)}$.

We prove first that $W(J''(\bs),k)\subseteq M^{(k)}_Q$.
Given $u$ in $W(J''(\bs),k)$, there is
$\bc\in\bZ^r$ such that $\sum_jc_j=k$, $c_j\geq 0$ for 
$j$ in $J''(\bs)$, and $u-\sum_jc_jv_j$ is in $\bZ_{>0}^n$.
If $v_{0}$ is in $\Gamma_{\fa}\cap Q$, we see that
$u-kv_0-e$ is in the semigroup generated by $v-w$,
where $v\in\Gamma_{\fa}$ and $w\in\Gamma_{\fa}\cap Q$. 
Therefore $u$ is in $M^{(k)}_Q$.

We show now the reverse inclusion.
It follows from definition that for every $j$ we have
$v_j+W(J''(\bs),k)\subseteq W(J''(\bs),k+1)$.
Therefore we may assume
$ k = 0 $, and it is enough to show that if $v$ is in $\Gamma_{\fa}$
and $w$ is in $\Gamma_{\fa}\cap Q$, then there is $\bc\in\bZ^r$
such that $\sum_jc_j=0$, for $j\in J''(\bs)$ we have
$c_j\geq 0$, and $v-w-\sum_jc_jv_j$ is in $\bZ_{\geq 0}^n$.
Since for every $v$ as above we can find $j\leq r$ such that
$v-v_j$ is in $\bZ_{\geq 0}^n$, we see that it is enough
to show the following: for every $w$ in $\Gamma_{\fa}\cap Q$,
there is $\bc\in\bZ^r$ such that $\sum_jc_j=-1$, for every $j\in J''(\bs)$
we have $c_j\geq 0$, and
$$-w\in\sum_{j=1}^rc_jv_j+\bZ_{\geq 0}^n.$$

By the choice of $Q$, we can find $\lambda_j\in\bQ_{\geq 0}$ for
$j\in J'(\bs)$ such that
$\sum_{j\in J'(\bs)}\lambda_j=1$ and $\sum_{j\in J'(\bs)}\lambda_jv_j$
is in the interior of $Q$. It follows that the convex cone generated
by 
$$\{v-\sum_{j\in J'(\bs)}\lambda_jv_j\mid v\in\Gamma_{\fa}\cap Q\}$$
is equal to the linear span of the same set.
On the other hand,
for every $v$ in $\Gamma_{\fa}\cap Q$, we write
$$v-\sum_{j\in J'(\bs)}\lambda_jv_j=\sum_{j\in J'(\bs)}\lambda_j(v-v_j),$$
so $v-\sum_{j\in J'(\bs)}\lambda_jv_j$ is in the convex cone generated by
the $v-v_j$, with $j\in J'(\bs)$. 
These two facts imply that given $w\in\Gamma_{\fa}\cap Q$, we can find
$m$ divisible enough
and for $j\in J'(\bs)$ elements $q_j$ in $\bZ_{\geq 0}$,
and $u_j$ in $\Gamma_{\fa}\cap Q$ such that the
$ m\lambda_j $ are integers and
$$m\left(\sum_{j\in J'(\bs)}\lambda_jv_j-w\right)=
\sum_{j\in J'(\bs)}q_j(u_j-v_j).$$
Therefore we can write
$$-w=(m-1)w+\sum_{j\in J'(\bs)}(q_ju_j-(q_j+m\lambda_j)v_j).$$
After replacing $w$ and each $u_j$ by suitable 
elements in $\{v_1,\ldots,v_r\}$, we get
$\bc\in\bZ^r$ such that $\sum_jc_j=-1$, $c_j\geq 0$ for $j\in J'(\bs)$
and
$-w\in\sum_jc_jv_j+\bZ_{\geq 0}^n$. This completes the proof of this step.

\smallskip

\noindent{\emph{Step 7}}. We show now that if $\gamma$ is a root
of $b_{\fa}$, then there is $Q$ such that $\gamma$ is in $R_Q$.
We have seen that there is $\bs\in\bQ^r$ such that
$-\bs$ is in $V(I_{\fa})$ and $\gamma=-\sum_js_j$. 
Moreover, by \emph{Step 3}
we may assume that $\sum_js_jv_j$ is in $\bZ_{>0}^n$.
If $Q$ is the smallest face of $P_{\fa}$ containing all
$v_j$ for $j\in J'(\bs)$, $Q$ is a proper face by \emph{Step 5}.

Since $\sum_{j=1}^rs_jv_j$ is in $\bZ_{>0}^n$, we see that
for every $i$ there is $j\in J'(\bs)$ such that $a_{i,j}>0$.
Therefore $Q$ is not contained in the hyperplane $(x_i=0)$,
for any $i$. 

Using \emph{Steps 4} and \emph{6}, we deduce that
$\sum_{j=1}^rs'_jv_j$ is in $M_Q^{(k(\bs))}\smallsetminus
M_Q^{(k(\bs)+1)}$. Moreover, this vector lies also 
in the linear space spanned by $Q$ (since $s'_j=0$ if $j\not\in J'(\bs)$),
so it follows from Remark~\ref{rem6} that
$-\sum_{j=1}^r s'_jL_Q(v_j)-k(\bs)$ is in $R_Q$.
If $j$ is in $J'(\bs)$, then $s'_j=s_j$ and $L_Q(v_j)=1$,
so $-\sum_{j\in J'(\bs)}s_j-\sum_{j\in J''(\bs)}s_j=\gamma$ is in $R_Q$.

\smallskip

\noindent{\emph{Step 8}}. We prove now the converse:
every element of some $R_Q$ is a root of $b_{\fa}$. 
By Remark~\ref{rem4} and the argument in \emph{Step 4}, we 
may assume that $Q$ is bounded. In particular, the convex cone spanned by
$Q$ is spanned by those $v_j$ in $Q$.

Suppose that $u$ is in $(M_Q\smallsetminus M_Q')\cap V_Q$.
We can write
$$u+\sum_{j=1}^ra_jv_j=\sum_{j=1}^rb_jv_j,$$
for some $a_j\in\bZ_{\geq 0}$
and $b_j\in\bQ_{\geq 0}$ such that $a_j=b_j=0$ if $v_j$ is not in $Q$
(we use the fact that given an element in $V_Q$, its sum
with a suitable positive integral combination of the $v_j$ in $Q$
lies in the convex cone generated by the $v_j$ in $Q$).
Put $q=\sum_ja_j$. 

Let $\bs\in\bQ^r$ be given by $s_j=b_j-a_j$ for all $i$. 
Since $b_j=0$ if $v_j$ is not in $Q$, it follows that the smallest face
$Q'$ of $P_{\fa}$ containing all  $v_j$ with $j\in J'(\bs)$ is contained in $Q$.
For every $j$, we can write $s'_j=b_j-p_j$ for some $p_j\in\bZ$.
Note that if $v_j$ is not in $Q$, then $p_j=0$.
Since
$ s''_j = a_j - p_j $,
we have $k(\bs)=\sum_ja_j-\sum_jp_j=q-p$, where $p=\sum_jp_j$.

By hypothesis $u$ is not in $M_Q'$ and
$$\sum_js'_jv_j=u+\sum_ja_jv_j-\sum_jp_jv_j\not\in M_Q^{(q-p+1)},$$
hence by \emph{Step 4} $-\bs$ is in $V(I_{\fa})$.
Since $L_Q=1$ on $Q$, we deduce $L_Q(u)+q=\sum_jb_j$.
Therefore $-L_Q(u)=\sum_ja_j-\sum_jb_j=-\sum_js_j$,
so $-L_Q(u)$ is a root of $b_{\fa}$. This completes the proof of this step
and that of the theorem.
\end{proof}

\begin{remark}\label{rem8}
For every facet $Q$ of $P_{\fa}$ that is not contained in any coordinate
hyperplane,
we get a root of $b_{\fa}$ given by $-L_Q(e)$.
Indeed, $e$ is clearly in $M_Q\cap V_Q$, so it is enough to show
that $e$ is not in $M_Q'$. Since $L_Q\geq 1$ on $P_{\fa}$, with equality
on $Q$, we see that $L_Q(u)\geq L_Q(e)+1$ for every $u$ in $M_Q'$.
This shows that $e$ is not in $M_Q'$. 
\end{remark}

\begin{remark}\label{negativity}
It is a general fact that the roots of the Bernstein-Sato polynomial
of an arbitrary ideal are negative rational numbers (see \cite{BMS2}).
In the case of a monomial ideal, this follows from Theorem~\ref{main1}.
Indeed, 
given a face $Q$ of $P_{\fa}$ that is not contained in any coordinate
hyperplane,
we can choose $L_Q$ such that $L_Q> 1$ on $P_{\fa}\smallsetminus Q$. 
Therefore for every $u$ in $M_Q$ we have $-L_Q(u)\leq -L_Q(e)$.
Moreover, if $m$ is large enough, then $me$ is in $P_{\fa}\smallsetminus Q$,
so $L_Q(e)>0$.
\end{remark}

\begin{remark}\label{rem9}
If $n\geq 2$, in Theorem~\ref{main1} it is enough to consider only those
faces $Q$ of $P_{\fa}$ of positive dimension. In order to show this,
we start with a more general remark: suppose that $Q\subset Q_1$
are faces of $P_{\fa}$ such that
$Q$ is not contained in any coordinate hyperplane.
Suppose that $u$ is in $(M_Q\smallsetminus M_Q')\cap V_Q$.
Then $u$ clearly lies in $M_{Q_1}\cap V_{Q_1}$. If $u$ is not in
$M_{Q_1}'$, then the roots corresponding to $u$ in $R_Q$ and $R_{Q_1}$
are the same. If this is not the case, then there is a root of $b_{\fa}$
in $R_{Q_1}$ of the form 
$-L_{Q}(u)+k$
for some $k\in\bZ_{>0}$. In particular, it follows from Remark~\ref{negativity}
that if $-L_{Q}(u)\geq -1$, then $-L_Q(u)$ is in $R_{Q_1}$.
However, if $Q$ is a vertex of $P_{\fa}$, then
$ u $ is in the segment between the origin and the vertex
and clearly
$ L_{Q}(u)\le 1 $.
\end{remark}

We give now the proof of the description of the classes
mod $\bZ$ of the roots of $b_{\fa}$.

\begin{proof}[Proof of Theorem~\ref{main2}]
We need to describe the set of classes mod $\bZ$ of the elements in $R_Q$,
for the faces $Q$ of $P_{\fa}$ that are not contained in coordinate
hyperplanes.
It follows from Remark~\ref{rem6} that 
this is equal to the set of classes
mod $\bZ$ of $-L_Q(w)$, where $w$ varies over $M_Q\cap V_Q$.
After replacing $Q$ by a facet containing it, it is enough to consider the
case when $Q$ is a facet, so $V_Q=\bR^n$.

Moreover, given $w\in \bZ^n$, there is a positive integer
$ k $ together with $u\in\Gamma_{\fa}\cap Q$ such that
$ w + ku $ is contained in
$ M_Q $. Since
$L_Q(w+ku)\equiv L_Q(w)$ (mod $\bZ$), we need to compute
the set of classes mod $\bZ$ of $-L_Q(w)$ when $w$ is in $\bZ^{n}$.
Write
$ L_Q(x)=\sum_i\beta_ix_i/m_Q $ where the
$ \beta_i $ are integers whose greatest common divisor is one.
Then the assertion is clear.
\end{proof}

We mention that unlike in Theorem~\ref{main2}, 
in the statement of Theorem~\ref{main1}
it is not enough to restrict to the facets of $P_{\fa}$
(see Example~\ref{non_maximal}
below). 

\begin{remark}\label{notation}
The formula in Theorem~\ref{main1} is not very convenient for
explicit computations since the set $(M_Q\smallsetminus M'_Q)\cap V_Q$
is infinite. We introduce now some finite sets that will be useful in the
next section when looking at examples.

For every face $Q$ that is not contained in any coordinate hyperplane
and for every $k\in\bZ_{\geq 0}$ we have
$$R_Q=\{-L_Q(u)+k\mid u\in (M_Q^{(k)}\smallsetminus M_Q^{(k+1)})\cap V_Q\}$$
Since $Q$ is not contained in a coordinate hyperplane, it follows that
for every $u$ in $V_Q\cap\bZ^n$, 
we can find $w_1,\ldots,w_s$ in $\Gamma_{\fa}\cap Q$
such that $u':=u+w_1\ldots+w_s$ is in $\bZ_{>0}^n$.
It is clear that $u$ is in $M_Q^{(k)}$ if and only if
$u'$ is in $M_Q^{(k+s)}$
and $-L_Q(u)+k=-L_Q(u')+k+s$.
Let us denote by $G_Q$ the subgroup of $\bZ^n$
generated by $\Gamma_{\fa}\cap Q$ and let
$G_Q':=w+G_Q$, where $w$ is in $\Gamma_{\fa}\cap Q$
(this clearly does not depend on $w$).
After subtracting a suitable number of elements of 
$\Gamma_{\fa}\cap Q$ from $u'$, we arrive at $u''$
that is in $\bZ_{>0}^n$, but not in $\bZ_{>0}^n+G_Q'$.

We want to find convenient subsets $E_Q$ of $V_Q$
such that if we put $E_Q^{(k)}:=E_Q\cap (M_Q^{(k)}\smallsetminus
M_Q^{(k+1)})$, we have
$$R_Q=\bigcup_{k\in\bZ_{\geq 0}}\{-L_Q(u)+k\mid u\in E_Q^{(k)}\}.$$
For example, the previous discussion shows that we may take
$E_Q=(\bZ_{>0}^n\cap V_Q)\smallsetminus (G_Q'+ \bZ_{> 0}^n)$.
One can check that if $Q$ is bounded, then $E_Q$ is finite.
It may be conjectured that $L_Q(u)\leq n$ for every $u\in E_Q$,
so $E_Q^{(k)}$ is empty for $k\geq n$ by Remark~\ref{negativity}.
This is easily proved if $n=2$. On the other hand, Example~\ref{optimal} below
shows that the conjecture is optimal: for every $n\geq 2$ there are examples
such that $E_Q^{(n-1)}$ is nonempty.

Sometimes we can choose $E_Q$ better: suppose for example
that we have a finite set
$\{p_i\}_i$ such that 
$$(V_Q\cap\bZ_{>0}^n)\smallsetminus (\bZ_{>0}^n+G_Q')
\subseteq\bigcup_i(p_i-(\bZ_{\geq 0}\cap V_Q)+G_Q).$$
In this case we may replace the above $E_Q$ by 
the union of the $E_Q\cap (p_i-(\bZ_{\geq 0}\cap V_Q))$.

Consider for example the case when $n=2$ 
(note that by Remark~\ref{rem9} it is enough the consider
only the one-dimensional faces $Q$). 
If $v=(a,b)$ and $v'=(a',b')$ are elements of
$Q\cap G_Q'$ such that $a<a'$, $b>b'$ and
$G_Q$ is generated by $v-v'$, then we may
replace the above $E_Q$ by $\{(i,j)\in\bZ^2_{>0}\mid i\leq a',j\leq b\}$.

Note that if $n=2$, then the roots corresponding to unbounded faces
of $P_{\fa}$ can be described as follows. Suppose, for example that 
$Q=\{(a,y)\mid y\geq b\}$ for integers $a$ and $b$, with $a>0$.
Then we have $R_Q=\{-k/a\mid 1\leq k\leq a\}$.
\end{remark}

\begin{remark}
(i) It is known that if we restrict to the
interval $(0,1)$,
the jumping coefficients of a monomial ideal
$ \fa $
coincide with those of a generic polynomial
$ f $ with the same
Newton polyhedron, see \cite{Laz} and \cite{How1}.
Up to a sign, these jumping coefficients are roots of the Bernstein-Sato polynomial
$ b_{f}(s) $ (see \cite{ELSV}) and also roots of
$ b_{\fa}(s) $ (see \cite{BMS2}).
However, it does not seem easy to describe in general
the relation between $b_{\fa}$ and $b_f$
(for a simple case, see the example discussed below).

\medskip
(ii) Consider the case of a weighted homogeneous polynomial
$ f $ with
weights
$ (w_{1},\dots,w_{n}) $ such that
$ a_{i}: = 1/w_{i} $ are positive integers.
If $f$ has an isolated singularity, then
the roots of
$ b_{f}(s)/(s + 1) $ are
$$
- \sum_{i = 1}^{n}\frac{p_{i}}{a_{i}}\quad\text{for}\quad
1\le p_{i}\le a_{i} - 1,
$$
with multiplicity one.
This is due to Kashiwara, and it also follows from a well-known
theorem of Malgrange
\cite{Mal} together with
a calculation of the Gauss-Manin connection due to Brieskorn.
We mention that a similar formula for the spectrum was given by
Steenbrink in \cite{Ste1}.
Note, however, that if we add monomials of higher degree to $f$, then
we might need to shift the above roots. For example, if 
$f=x^5+y^4$ and $g=x^5+y^4+x^3y^2$, then the root 
$-\frac{31}{20}$ of $b_f$ is shifted to give the root $-\frac{11}{20}$
of $b_g$.

On the other hand, if $\fa=(x_1^{a_1},\ldots,x_n^{a_n})$, then
the roots of
$ b_{\fa}(s) $ are
$$
- \sum_{i = 1}^{n}\frac{p_{i}}{a_{i}}\quad\text{for}\quad
1\le p_{i}\le a_{i}.
$$
However, if
$ \fa $ has extra generators inside or on the boundary of
the Newton polyhedron, then 
we might need to shift the above roots
(see Example~\ref{ex2} below).

\medskip
(iii) In the case
$ n = 2 $ there is a similar formula
due to Steenbrink \cite{Ste2} for the spectrum of
a generic function having the Newton polygon
$ P_{\fa} $ whose complement in
$ \bR_{\ge 0}^{n} $ is bounded (see also \cite{Sai} for the case
$ n > 2 $).
His formula is in terms of
$\beta(u):=\min_QL_Q(u)$, where $Q$ varies over the compact facets
of the Newton polygon.
The spectral numbers that are $\leq 1$ are given by
$\beta(u)$ for those $u\in\bZ_{>0}^2$ with $\beta(u)\leq 1$.
The remaining spectral numbers are obtained by symmetry:
if $\alpha_1\leq\ldots\leq\alpha_{\mu}$ are all the spectral numbers, then
$\alpha_i+\alpha_{\mu+1-i}=2$.

Restricting to the numbers
less than $1$, Steenbrink's formula is the same as the
one for the jumping coefficients 
obtained from \cite{How1} and \cite{How2}.
Note that if we further restrict to those $u$ in the
cone over the facet $Q$ of $P_{\fa}$, then only 
the
integral points inside this cone are taken into account.

However, recall that if $v=(a,b)$ and $v'=(a',b')$ are points
on $Q$ as in
Remark~\ref{notation}, then 
our formula considers all the integral points in
$$
E_{Q_i} = [1,a']\times[1,b].
$$
For example,
$ L_{Q}(e) $ is always roots of the
Bernstein-Sato polynomial,
but usually it is neither a jumping coefficients nor a spectral number
(unless
$ \partial P_{\fa} $ has only one compact face or it has two such faces
and
$ e $ is on the middle one-dimensional cone).
\end{remark}

\section{Examples}

In this section we give some examples to 
illustrate our combinatorial description.
We freely use the notation introduced in Remark~\ref{notation}.

\begin{example}\label{ex1}
Let $\fa=(x^ay,xy^b)$, with
$ a,b \ge 2 $.
We see that the roots are
$$
- \frac{(b - 1)i + (a - 1)j}{ab - 1}\quad
\text{for}\quad 1\le i\le a,\,\, 1\le j\le b,
$$
and we include also $-1$ if it does not appear
in this list.
Indeed, if $Q$ is the only bounded one-dimensional face of $P_{\fa}$,
then we may take
$$
E_Q= E_{Q}^{(0)} = \{(i,j)\in\bZ_{>0}^{2}\,|\,
i \le a,j \le b\},
$$
(see Remark~\ref{notation}),
and
$ L_{Q}(x,y) = ((b - 1)x + (a - 1)y)/(ab - 1) $.
\end{example}

\begin{example}\label{ex2}
Let $\fa=(x^d,x^{d-a}y^a,y^d)$, where 
$a$ and $d$ are positive integers such that $d\geq 2$ and
$d/a\in\bZ$. If $Q$ is the only bounded one-dimensional face
of $P_{\fa}$, then we may take
$$
E_{Q} = E_{Q}^{(0)} = \{(i,j)\in\bZ_{>0}^{2}\,|\,
i \le d,j \le a\}.
$$
Since
$ L_Q(x,y) = (x + y)/d $, we see that the roots of $b_{\fa}$ are
$ \{ - k/d\,|\,2\le k\le d + a\} $.
This illustrates how the roots depend on
$ \Gamma_{\fa}\cap Q$.
\end{example}

\begin{example}\label{ex3}
If $\fa=(xy^5,x^3y^2,x^4y)$, then 
the roots of $b_{\fa}$ are
$$
-\frac{i}{13}\,\,(5\le i\le 17),\,\,-\frac{j}{5}\,\,(2\le j\le 6).
$$
This is compatible with a calculation by
\emph{Macaulay}2 using the method in \cite{BMS2}
(in this case the ideal
$ I_{\fa} $ introduced in \S 2 has more than 20 generators).
This is one of the simplest examples such that
$ E_{Q}^{(1)}\ne
\emptyset $.
Indeed, if $Q$ is the face 
of $P_{\fa}$ containing $(1,5)$ and $(3,2)$
then by Remark~\ref{notation} we may take
$$
E_{Q} = \{(i,j)\in\bZ_{>0}^{2}\,|\,i\le 3,j\le 5\},
$$
and
$ L_{Q}(x,y)= (3x + 2y)/13 $.
Therefore we get
$$
E_{Q}^{(1)} = \{(3,5)\}\ne\emptyset.
$$
and 
$-19/13 $ is shifted to
give the root$-6/13 $.
This can be compared with the next example.
\end{example}

\begin{example}\label{ex4}
Let $\fa$ be the ideal $(xy^5,x^3y^2,x^5y)$.
The roots of $b_{\fa}$ are
$$
-\frac{5}{13},\,\,-\frac{i}{13}\,\,(7\le i\le 17),\,\,
-\frac{19}{13},\,\,-\frac{j}{6}\,\,(3\le j\le 9).
$$
If $Q$ is the same face as in the previous example, 
$E_Q$ and $L_Q$ are the same as before, but $E_Q^{(k)}=\emptyset$
for all $k\geq 1$.
One can see that
$ R_{Q} $ depends also on the intersection of
$ \Gamma_{\fa} $ with
the faces of $P_{\fa}$ different from $Q$.
\end{example}

\begin{example}\label{non_maximal}
Let 
$ v_{i}\,(1\le i\le 4) $ be
$ (3,0,3) $,
$ (0,3,3) $,
$ (0,0,7) $,
$ (1,1,6)$
respectively, and consider the ideal $\fa$ generated by the $x^{v_i}$
in $\bC[x_1,x_2,x_3]$.
The roots of $b_{\fa}$ are
$$
-\frac{7}{21},\,-\frac{11}{21},
\,\,-\frac{i}{21}\,\,(14\le i\le 39),\,\,-\frac{42}{21},
$$
where
$ -7/21 $ comes from the face of $P_{\fa}$ determined by
the cone generated by
$v_{1}$, $v_{2}$, $e_{1}$, and $e_{2} $.
Here
$ e_{i} $ denotes the
$ i $-th unit vector.
Let
$ Q_{1} $,
$ Q_{2} $ be the
faces determined by
$\{v_{1},v_{2},v_{3}\}$ and
$\{v_{1},v_{2}\}$, respectively.
By Remark~\ref{notation} we may take
$$
\aligned E_{Q_1}& = \{(i,j,k)\in\bZ_{>0}^{3}\,|\,
i,j\le 3,\,k\le 7\},\\
E_{Q_2}& = \{(i,j,k)\in\bZ_{>0}^{3}\,|\,
i,j\le 3,\,\,k = i + j,\},
\endaligned
$$
and
$ L_{Q_{1}}(x,y,z) = (4x + 4y + 3z)/21 $.
We have
$ E_{Q_{2}}^{(k)} = \emptyset $ for $k\geq 1$ and
$ E_{Q_{1}}^{(k)} = \emptyset $ for $k\geq 2$, and
$$
E_{Q_{1}}^{(1)} = \{(i,j,k)\in
E_{Q_{1}}\,|\,i,j\ge 2,\,k = 7\,\,
\text{or}\,\,i,j = 3,\,\,k = 6\}.
$$
Note that
$ u = (3,3,6) $ belongs to
$ E_{Q_{1}}^{(1)} $ and
$ E_{Q_{2}}^{(0)} $, and
$ L_{Q_{1}}(u) = 2 $.
Therefore
$ R_{Q_{2}} $ is not contained in
$ R_{Q_{1}} $, and we see that in Theorem~\ref{main1}
we need to consider also 
$ R_{Q} $ for faces $Q$ with $\dim(Q)<n-1$.
\end{example}

\begin{example}\label{ex6}
Let
$ v_{i}\,(1\le i\le 4) $ be
$ (4,0,1) $,
$ (0,3,2) $,
$ (0,0,5) $,
$ (1,1,3)$, respectively, and let $\fa$ be generated by the $x^{v_i}$
in $\bC[x_1,x_2,x_3]$.
Let
$ Q_{1} $,
$ Q_{2} $,
$ Q_{3} $ be the faces determined by
the cones generated by
$ (v_{1},v_{2},v_{3}) $,
$ (v_{1},v_{2},e_{2}) $ and
$ (v_{1},v_{2}) $, respectively.
In this case, all
$ v_{i} $ lie on
$ Q_{1}$ and
$ G_{Q_{1}}^{(1)} = \bZ^{3}\cap L_{Q_{1}}^{-1}(1) $.
Therefore we may take
$$
\aligned E_{Q_{1}}& = \{(i,j,k)\in\bZ_{>0}^{3}\,|\,
i + j + k\le 7\},\\
E_{Q_{2}}& = \{(i,j,k)\in\bZ_{>0}^{3}\,|\,
i\le 4,\,\,k\le 2\},\\
E_{Q_{3}}& = \{(4,3,3)\},\endaligned
$$
and
$ L_{Q_{1}}(x,y,z) = (x + y + z)/5 $,
$ L_{Q_2}(x,y,z) = x/8 + z/2 $.
We have
$ E_{Q_{i}}^{(k)} = \emptyset $
for every $i$ and for every $k\geq 1$, so
$ R_{Q_{1}} $,
$ R_{Q_{2}} $ and
$ R_{Q_{3}} $ respectively consist
of
$$
- \frac{i}{5}\,\,(3\le i\le 7),\,\, - \frac{j}{8}\,\,(5\le j\le 12),\,\,
\text{and}\,\,\, - 2,
$$
and they give all the roots.
\end{example}

\begin{example}\label{optimal}
Let
$ v_{i}\,(1\le i\le 4) $ be
$ (3,0,0) $,
$ (0,3,0) $,
$ (0,0,5) $,
$ (1,1,2)$, and $\fa$ the corresponding monomial ideal.
If 
$ Q $ is the face of $P_{\fa}$ with vertices
$ (v_{1},v_{2},v_{3}) $, then we may take
$$
E_{Q} = \{(i,j,k)\in\bZ_{>0}^{3}\,|\,
i,j\le 3,\,k\le 5\},
$$
and
$ L_{Q}(x,y,z)= (5x + 5y + 3z)/15 $.
In this case we have
$$
\aligned
E_{Q}^{(2)}& =\{(3,3,5)\}\ne\emptyset,\,\,\,\,E_Q^{(k)}=0\,{\rm for}\,k\geq 3,\\
E_{Q}^{(1)}\cup E_Q^{(2)} & = \{(i,j,k)\in\bZ_{>0}^{3}\,|\,i,j\ge 2,\,k\ge 3\}\cup
\{(3,1,5),(1,3,5)\},
\endaligned
$$
and the roots are
$$
-\frac{i}{15}\,(13\le i\le 32).
$$
Note that
$ (3,1,5) \in E_{Q}^{(1)} $ comes from
$ (3,0,5) = v_{2} + 2(v_{4} - v_{2}) + e $, and
$ -35/15 $ is shifted to
$ -20/15 $.

This example can be generalized to any
$ n \ge 3 $ by taking $v_i=ne_i$ for $1\leq i\leq n-1$,
$v_n=(2n-1)e_n$, and $v_{n+1}=(1,\ldots,1,2)$.
\end{example}

%\newpage

\end{document}